\documentclass[10pt]{llncs}
\usepackage{graphics,amsmath,amssymb}
\def\Q{{\bf Q}}
\def\Z{{\bf Z}}
\def\Sym{{\rm Sym}}
\def\polhk{\c}
\DeclareFontEncoding{OT2}{}{} 
  \newcommand{\textcyr}[1]{%
    {\fontencoding{OT2}\fontfamily{wncyr}\fontseries{m}\fontshape{n}%
     \selectfont #1}}

\begin{document}
\title{Symmetric powers of elliptic curve L-functions}
\author{
Phil Martin and
Mark Watkins\inst{}\thanks{Supported during parts of this research
  by EPSRC grants GR/N09176/01 and GR/T00658/01,
  the Isaac Newton Institute, the CNRS and the Institut Henri Poincar\'e,
  and the MAGMA Computer Algebra Group at the University of Sydney.}}
\institute{University of Bristol\\
\email{phil\_martin\_uk@hotmail.com\ \ watkins@maths.usyd.edu.au}}

\maketitle
\begin{abstract}
The conjectures of Deligne, Be\u\i linson, and Bloch-Kato assert that
there should be relations between the arithmetic of algebro-geometric
objects and the special values of their $L$-functions.
We make a numerical study for symmetric power $L$-functions
of elliptic curves, obtaining data about the validity
of their functional equations, frequency of vanishing of
central values, and divisibility of Bloch-Kato quotients.
\end{abstract}

\section{Introduction and Motivation}
There are many conjectures that relate special values of $L$-functions
to the arithmetic of algebro-geometric objects.
The celebrated result $\zeta(2)=\pi^2/6$ of Euler~\hbox{\cite[\S XV]{euler}}
can be reinterpreted as such, but Dirichlet's class
number formula~\hbox{\cite[\S 5]{dirichlet}} is better seen
to be the primordial example.
Modern examples run the gamut, from conjectures of Stark \cite{stark}
on Artin $L$-functions and class field theory,
to that of Birch and Swinnerton-Dyer \cite{BSD} for elliptic curves,
to those of Be\u\i linson \cite{b1,beilinson} related to $K$-theory,
with a passel of others we do not mention.
For maximal generality the language of motives
is usually used (see \cite[\S 1-4]{flach}).

One key consideration is where the special value is taken.
The $L$-function can only vanish inside the critical strip
or at trivial zeros; indeed, central values (at the center
of symmetry of the functional equation) are the most interesting ones
that can vanish, and the order of vanishing is likely
related to the rank of a geometric object
(note that orders of trivial zeros can be similarly interpreted).

We have chosen to explore a specific family of examples, namely
symmetric power $L$-functions for rational elliptic curves.
The impetus for this work was largely a theoretical result
\cite{martin} of the first author,
whose computation of Euler factors in the difficult case of additive primes
greatly reduced the amount of hassle needed to do large-scale computations.
Previous theoretical work includes that of Coates and Schmidt
\cite{coates-schmidt} on the symmetric square and Buhler, Schoen, and Top
on the symmetric cube~\cite{BST}; this second paper also contains
a lot of computational evidence, while Watkins has provided
much data \cite{watkins} in the symmetric square case.
In some cases the \hbox{$L$-functions} we use
are not known to possess the properties that would be required
to justify that our computations produce numbers of any
validity whatsoever --- in these cases, the ``numerical coincidence''
in our computations can be seen as evidence for the relevant conjectures.

\subsection{Acknowledgements}
We thank Neil Dummigan for useful comments and
his interest in these computations;
indeed much of this paper came out of a desire to generalise
the numerical experiment in~\hbox{\cite[\S 7]{dummigan}.}
We thank Frank Calegari for pointing out that we have a
factorisation in the CM case, and Erez Lapid for help with local root numbers.

\section{$L$-functions}
We define the symmetric power $L$-functions of an elliptic curve~$E/\Q$
via computing an Euler factor at every prime~$p$. This Euler factor
is computed by a process that essentially just takes the symmetric
power representation of the standard \hbox{2-dimensional} Galois representation
associated to~$E$, and thus our method is a generalisation of that used
by Coates and Schmidt~\cite{coates-schmidt} for the symmetric square,\footnote{
 Note that Buhler, Schoen, and Top \cite{BST} phrase
 their definition of Euler factors differently,
 as they emphasise that conjecturally the $L$-function
 is related to a motive or higher-dimensional variety;
 however, their definition is really the same as ours.}
following the original description of Serre~\cite{serre:eulerproduct}.
We briefly review the theoretical framework, and give explicit
formulae for the Euler factors in a later section.

For every prime~$p$ choose an auxiliary prime~$l\neq 2,p$
and fix an embedding of $\Q_l$ into~${\bf C}$.
Let $E_t$ denote the $t$-torsion of~$E$,
and $T_l(E)=\lim\limits_\leftarrow E_{l^n}$
be the $l$-adic Tate module of~$E$ (we fix a basis).
The module $V_l(E)=T_l(E)\otimes_{\Z_l}\Q_l$
has dimension~2 over~$\Q_l$ and has a natural action
of~${\rm Gal}(\overline{\Q}_p/\Q_p)$
[indeed one of of~${\rm Gal}(\overline{\Q}/\Q)$],
and from this we get a representation
$\rho_l:{\rm Gal}(\overline{\Q}_p/\Q_p)\rightarrow {\rm Aut}(V_l)$.
We write $H_l^1(E)={\rm Hom}_{\Q_l}(V_l(E),\Q_l)$,
and take the $m$th symmetric power
of the contragredient of~$\rho_l$, getting
$$\rho_l^m:{\rm Gal}(\overline{\Q}_p/\Q_p)
\rightarrow {\rm Aut}\bigl(\Sym^m\bigl(H_l^1(E)\bigr)\bigr)
\subset GL_{m+1}({\bf C}).$$
We write $D_p={\rm Gal}(\overline{\Q_p}/\Q_p)$,
let $I_p$ be the inertia group of this extension,
and let ${\rm Frob}_p$ be the element of
$D_p/I_p\cong {\rm Gal}(\overline{{\bf F}}_p/{\bf F}_p)$
given by $x\rightarrow x^p$. With all of this, we have
$$L(\Sym^m E,s)=
\prod_p {\rm det}\Bigl[\text{Id}_{m+1}-\rho_l^m({\rm Frob}_p^{-1}) p^{-s}\Bigm|
\bigl(\Sym^m\bigl(H_l^1(E)\bigr)\bigr)^{I_p}\Bigr]^{-1}.$$
For brevity, we write $L_m(E,s)=L(\Sym^m E,s)$, and denote the
factors on the right side by~$U_m(p;s)$.
As mentioned by Coates and Schmidt \cite[p.~106]{coates-schmidt},
it can be shown that $U_m(p;s)$ is independent of our choices.
The analytic theory and conjectures concerning
these symmetric power $L$-functions are described in~\cite{shahidi}.
In particular, the above Euler product converges in a half-plane,
and is conjectured to have a meromorphic continuation to the whole
complex plane.

We also need the conductor $N_m$
of this symmetric power representation.
We have $N_m=\prod_p p^{f_m(p)}$ where $f_m(p)=\epsilon_m(I_p)+\delta_m(p)$.
Here $\epsilon_m(I_p)$ is the codimension of
$\bigl(\Sym^m\bigl(H_l^1(E)\bigr)\bigr)^{I_p}$
in~$\Sym^m\bigl(H_l^1(E)\bigr)$;
we shall see that it can be computed via a character-theoretic argument.
The wild conductor $\delta_m(p)$ is 0 unless $p=2,3$, when it
can be computed as in \cite[\S 2.1]{serre:eulerproduct}
or the appendix of \cite{coates-schmidt}.

\subsection{Critical values}
The work of Deligne~\cite[Prop.~7.7ff]{deligne} tells us
when and where to expect critical values; these are a subset
of the more-general special values, and are the easiest to consider.\footnote{
 Critical values conjecturally only depend
 on periods (which are local objects),
 while the more-general special values can also
 depend on (global) regulators from $K$-theory.}
When $m=2v$ with $v$ odd there is a critical value~$L_m(E,v+1)$
at the edge of the critical strip,
and when $m=2u-1$ is odd there is a critical central value~$L_m(E,u)$.
We let $\Omega_{+},\Omega_{-}$ be the real/imaginary periods of~$E$
for $m\equiv 1,2\pmod{4}$, and vice-versa for~$m\equiv 3\pmod{4}$.
In the respective cases of $m$ even/odd we expect rationality
(likely with small denominator) of either
\begin{equation}\label{eqn:special}
{L_m(E,v+1)\over (2\pi)^{v+1}}
\biggl({2\pi N\over\Omega_{+}\Omega_{-}}\biggr)^{v(v+1)/2}
\qquad\text{or}\qquad
{L_m(E,u)(2\pi N)^{u(u-1)/2}
\over\Omega_{+}^{u(u+1)/2}\Omega_{-}^{u(u-1)/2}}.
\end{equation}
When $m$ is odd, the order of $L_m(E,s)$ at $s=u$
should equal the rank of an associated geometric object.
The Bloch-Kato conjecture \cite{block-kato} relates the quotients
in \eqref{eqn:special} to $H^0$-groups, Tamagawa numbers, and generalised
Shafarevich-Tate groups.\footnote{
 See \hbox{\cite[\S 7]{dummigan}} for an explicit example;
 note his imaginary period is twice that of our normalisation
 (and the formula is out by a power-of-2 in any case),
 and the conductor enters the formula in a different place
 (this doesn't matter for semistable curves).}

\section{Computation of Euler factors and local conductors}
We first consider multiplicative and potentially multiplicative reduction
for a given prime~$p$; these cases can easily be detected since $v_p(j_E)$,
the valuation of the $j$-invariant, is negative,
with the reduction being potentially multiplicative when~$p|c_4$.
When $E$ has multiplicative reduction,
the filtration of~\hbox{\cite[\S 8]{BST}} implies
the local tame conductor $\epsilon_m$ is~$m$
and $\delta_m(2)=\delta_m(3)=0$ for all~$m$.
The Euler factor is $U_m(p;s)=(1-a_p^m/p^s)^{-1}$,
where $a_p=\pm 1$ is the trace of Frobenius.
In the case of potentially multiplicative reduction,
for $m$ odd we have $\epsilon_m=m+1$, and so~$U_m(p;s)\equiv 1$,
while with $m$ even, we have that $\epsilon_m=m$ and
compute that $U_m(p;s)=(1-1/p^s)^{-1}$.
The wild conductor at $p=2$ is $\delta_m(2)={m+1\over 2}\delta_1(2)$
for odd~$m$ and is zero for even~$m$, while $\delta_m(3)=0$ for all~$m$.

\subsection{Good and additive reduction --- tame conductors}
Let $E$ have good or potentially good reduction at a prime~$p$,
and choose an auxiliary prime $l\neq 2,p$. The inertia group $I_p$
acts on $V_l(E)$ by a finite quotient in this case.
Let ${\bf G}_p={\rm Gal}\bigl(\Q_p(E_l)/\Q_p)$ and $\Phi_p$
be the inertia group of this extension.\footnote{
 The group $\Phi_p$ is independent of the choice of~$l$
 (see~\cite[p.~312]{serre}), while only whether ${\bf G}_p$
 is abelian matters, and this independence follows
 as in \cite[Lemmata~1.4 \& 1.5]{coates-schmidt}.}
The work of Serre~\cite{serre} lists the possibilities for~$\Phi_p$.
It can be a cyclic group $C_d$ with $d=1,2,3,4,6$;
additionally, when $p=2$ it can be $Q_8$ or $SL_2({\bf F}_3)$,
and when $p=3$ it can be $C_3\rtimes C_4$.
For each group there is a unique faithful \hbox{2-dimensional}
representation $\Psi_\Phi$ of determinant~1 over~${\bf C}$,
which determines~$\overline\rho_l$.

Our result now only depends on~$\Phi$;
for a representation $\Psi$ we have the trace relation
(which is related to Chebyshev polynomials of the second kind)
\vspace{-1ex}
\begin{equation}
{\rm tr}(\Sym^m\Psi)=\sum_{k=0}^{m/2}
{m-k\choose k}{\rm tr}(\Psi)^{m-2k}(-{\rm det}\,\Psi)^k,\label{eqn:sigma}
\end{equation}
and from taking the inner product of ${\rm tr}(\Sym^m\Psi_\Phi)$
with the trivial character we find the dimension of the $\Phi$-fixed subspace
of~$\Sym^m\bigl(H_l^1(E)\bigr)$, which we denote by~$\beta_m(\Phi)$.
Upon carrying out this calculation,
we obtain Table~\ref{tbl:betaPhi}, which lists values for~$\beta_m(\Phi)$,
from which we get the tame conductor
$\epsilon_m(\Phi)=m+1-\beta_m(\Phi)$.
The wild conductors $\delta_m(p)$ are 0 for $p\ge 5$,
and for $p=2,3$ are described below.

\begin{table}[h]
\caption{Values of $\beta_m(\Phi)$ for various inertia groups;
here $\tilde m$ is $m$ modulo~12.\label{tbl:betaPhi}}
\begin{center}
\begin{tabular}{*{8}{|@{\hskip3.7pt}c@{\hskip3.7pt}}|}\hline
$\tilde m$&$C_2$&$C_3$&$C_4$&$C_6$&$Q_8$&
$C_3\rtimes C_4$&$SL_2({\bf F}_3)$\\\hline

$0$&$m+1$&$(m+3)/3$&$(m+2)/2$&$(m+3)/3$&$(m+4)/4$&$(m+6)/6$&$(m+12)/12$\\
$1$&$0$&$(m-1)/3$&$0$&$0$&$0$&$0$&$0$\\
$2$&$m+1$&$(m+1)/3$&$m/2$&$(m+1)/3$&$(m-2)/4$&$(m-2)/6$&$(m-2)/12$\\
$3$&$0$&$(m+3)/3$&$0$&$0$&$0$&$0$&$0$\\
$4$&$m+1$&$(m-1)/3$&$(m+2)/2$&$(m-1)/3$&$(m+4)/4$&$(m+2)/6$&$(m-4)/12$\\
$5$&$0$&$(m+1)/3$&$0$&$0$&$0$&$0$&$0$\\
$6$&$m+1$&$(m+3)/3$&$m/2$&$(m+3)/3$&$(m-2)/4$&$m/6$&$(m+6)/12$\\
$7$&$0$&$(m-1)/3$&$0$&$0$&$0$&$0$&$0$\\
$8$&$m+1$&$(m+1)/3$&$(m+2)/2$&$(m+1)/3$&$(m+4)/4$&$(m+4)/6$&$(m+4)/12$\\
$9$&$0$&$(m+3)/3$&$0$&$0$&$0$&$0$&$0$\\
$10$&$m+1$&$(m-1)/3$&$m/2$&$(m-1)/3$&$(m-2)/4$&$(m-4)/6$&$(m-10)/12$\\
$11$&$0$&$(m+1)/3$&$0$&$0$&$0$&$0$&$0$\\\hline
\end{tabular}
\end{center}
\vspace*{-8ex}
\end{table}

\subsection{Good and additive reduction --- Euler factors for~$p\ge 5$}
When $p\ge 5$, a result of Serre~\cite{serre}
tells us that the inertia group is~$\Phi=C_d$
where $d=12/\gcd\bigl(12,v_p(\Delta_E)\bigr)$.
Note that this gives $d=1$ when $p$ is a prime of good reduction,
which we naturally include in the results of this part.
We summarise the results of Martin's work \cite{martin}
concerning the Euler factors. Note that the result of
\hbox{D\polhk abrowski}~\hbox{\cite[Lemma~1.2.3]{dabrowski}}
appears to be erroneous.

There are two different cases for the behaviour of the Euler factor,
depending on whether the decomposition group
${\bf G}_p={\rm Gal}\bigl(\Q_p(E_l)/\Q_p\bigr)$ is abelian.
From~\cite[Prop.~2.2]{rohrlich} or~\cite[Th.~2.1]{watkins},
we get that this decomposition group is abelian
precisely when~$p\equiv 1\pmod{d}$.
When ${\bf G}_p$ is nonabelian we have
\begin{equation}
U_m(p;s)=(1-(-p)^{m/2}/p^s)^{-A_m}(1+(-p)^{m/2}/p^s)^{-B_m},
\label{eqn:Upnoncyclic}
\end{equation}
where $A_m+B_m=\beta_m$ and $A_m$ is the dimension of
$\bigl(\Sym^m\bigl(H_l^1(E)\bigr)\bigr)^{{\bf G}_p}$.
Using ${\bf G}_p/\Phi_p\cong C_2$ and ${\rm det}\bigl(\Psi_\Phi(x)\bigr)=-1$
for~$x\in{\bf G}_p\backslash\Phi_p$, more character calculations tell us
this dimension is $(\beta_m+1)/2$ when $\beta_m$ is odd
and is $\beta_m/2$ when $\beta_m$ is even.
This also holds for the non-cyclic~$\Phi$ when~$p=2,3$,
for which~${\bf G}_p$ is automatically nonabelian.
When $\Phi=C_3$ and $m$ is odd,
we have $U_m(p;s)=(1+p^m/p^{2s})^{-\beta_m/2}$.

When ${\bf G}_p$ is abelian,
we need to compute a Frobenius eigenvalue~$\alpha_p$
(whose existence follows from \cite[p.~499]{serre-tate}).
In the case of good reduction, this comes
from counting points mod~$p$ on the elliptic curve;
we have $\alpha_p=(a_p/2)\pm i\sqrt{p-a_p^2/4}$
where $p+1-a_p$ is the number of (projective) points on $E$ modulo~$p$.
And when $\Phi=C_2$ we count points on the $p$th quadratic twist of~$E$.
In general, we need to re-scale the coefficients of our curve by some power
of~$p$ that depends on the valuations $v_p$ of the coefficients.
Since $p\ge 5$, we can write our curve as \hbox{$y^2=x^3+Ax+B$},
and then re-scale by a factor $t=p^{\min(v_p(A)/2,v_p(B)/3)}$
to get a new curve \hbox{$E^t:y^2=x^3+Ax/t^2+B/t^3$,}
possibly defined over some larger field.
Because of our choice of~$t$, at least one of $A/t^2$ and $B/t^3$
will have $v_p$ equal to~0.
The reduction~$\tilde E^t$ modulo some (fractional) power of~$p$
is then well-defined and non-singular, and we get $\alpha_p$ from
counting points on~$\tilde E^t$; it turns out that choices of roots
of unity will not matter when we take various symmetric powers.
Returning back to $U_m(p;s)$, we get that when
$\Q_p(E_l)/\Q_p$ is abelian this Euler factor is
\begin{equation}
U_m(p;s)=
\prod_{0\le i\le m\atop d|(2i-m)}(1-\alpha_p^{m-i}\bar\alpha_p^i/p^s)^{-1}.
\label{eqn:Upcyclic}
\end{equation}

\subsection{Considerations when $p=3$}
Next we consider good and additive reduction for~$p=3$.
We first determine the inertia group,
using the 3-valuation of the conductor as our main guide.
In the case that $v_3(N)=0$ we have good reduction,
while when $v_3(N)=2$ and $v_3(\Delta)$ is even we have $\Phi=C_2$.
Since ${\bf G}_3$ is abelian here,
the Euler factor is given by~\eqref{eqn:Upcyclic},
while the wild conductor is~0 and tame conductor is obtained
from Table~\ref{tbl:betaPhi}.
When $v_3(N)=2$ and $v_3(\Delta)$ is odd we have that~$\Phi=C_4$
and ${\bf G}_3$ is nonabelian. The wild conductor $\delta_m(3)$ is~0,
and the Euler factor is given by~\eqref{eqn:Upnoncyclic}.

When $v_3(N)=4$ we get $\Phi=C_3$ or~$C_6$, the former
case when~$4|v_3(\Delta)$. For these inertia groups,
the question of whether ${\bf G}_3$ is abelian can
be resolved as follows (see~\cite[Th.~2.4]{watkins}).
Let $\hat c_4$ and $\hat c_6$
be the invariants of the minimal twist of~$E$ at~$3$.
In the case that $\hat c_4\equiv 9\pmod{27}$,
we have that ${\bf G}_3$ is abelian when $\hat c_6\equiv \pm 108\pmod{243}$
while if $3^3|\hat c_4$ then ${\bf G}_3$ is abelian
when~$\hat c_4\equiv 27\pmod{81}$.
In the abelian case we have $\alpha_3=\zeta_{12}\sqrt 3$ up to sixth roots,
which is sufficient.
The Euler factor is then given by either~\eqref{eqn:Upnoncyclic}
or~\eqref{eqn:Upcyclic},
the tame conductor can be obtained from Table~\ref{tbl:betaPhi},
and the wild conductor (computed as in the appendix of~\cite{coates-schmidt})
from Table~\ref{tbl:wild3}.
When $v_3(N)=3,5$ we have that~$\Phi=C_3\rtimes C_4$.
The Euler factor is given by~\eqref{eqn:Upnoncyclic}
and the wild conductor can be obtained from Table~\ref{tbl:wild3},
with the first $C_3\rtimes C_4$ corresponding to $v_3(N)=3$,
and the second to~$v_3(N)=5$.

\subsection{Considerations when $p=2$}
Finally we consider~$p=2$, where first we determine the inertia group.
Let $M$ be the conductor of the minimal twist $F$ of $E$ at~$2$, recalling
\hbox{\cite[\S~2.1]{watkins}} that in general we need
to check four curves to determine this twist.
Table~\ref{tbl:phi2} then gives the inertia group.
The appendix of \cite{coates-schmidt}
omits a few of these cases; see~\cite{watkins}.
When $\Phi=C_1,C_2$ we can always determine $\alpha_p$
via counting points modulo~$p$ on $E$ or a quadratic twist,
and ${\bf G}_2$ is always abelian.
The Euler factor is then as in~\eqref{eqn:Upcyclic}.
For $\Phi=C_3,C_6$ the group ${\bf G}_2$ is always nonabelian,
and the Euler factor is as in~\eqref{eqn:Upnoncyclic}.
For the case of $\Phi=C_4$ and $p=2$, the question
of whether ${\bf G}_2$ is abelian comes down
\hbox{\cite[Th.~2.3]{watkins}} to whether the $c_4$ invariant of $F$
is 32 or 96 modulo~128, it being abelian in the latter case,
where we have $\alpha_2=\zeta_8\sqrt 2$ up to fourth roots.
The Euler factors for this and the two cases of noncyclic~$\Phi$
are obtained from~\eqref{eqn:Upnoncyclic} or~\eqref{eqn:Upcyclic},
while the wild conductors $\delta_m(2)$
are given in Table~\ref{tbl:wild2},
with the appropriate line being determinable from
the conductor of the first symmetric power.

\vspace{-2.5ex}
\begin{table}[h]
\begin{minipage}{2.5in}
\caption
{Values for $\delta_m(2)$.\label{tbl:wild2}}
\begin{center}
\vspace{-2ex}
\begin{tabular}
{|@{\hskip4pt}c@{\hskip4pt}|*{1}{@{\hskip3.7pt}c@{\hskip3.7pt}}|
@{\hskip4pt}c@{\hskip4pt}|}\hline
$\Phi_2$&$m=1$&formula\\\hline
$C_2,C_6$&2&$\epsilon_m(C_2)$\\[0.5pt]
$C_2,C_6$&4&$2\epsilon_m(C_2)$\\[0.5pt]
$C_4$&6&$2\epsilon_m(C_4)+\epsilon_m(C_2)$\\[0.5pt]
$Q_8$&3&$\epsilon_m(Q_8)+{1\over 2}\epsilon_m(C_2)$\\[0.5pt]
$Q_8$&4&$\epsilon_m(Q_8)+\epsilon_m(C_2)$\\[0.5pt]
$Q_8$&6&
$\scriptstyle{{\epsilon_m(Q_8)+\epsilon_m(C_4)+\epsilon_m(C_2)}}$\\[1.5pt]
$SL_2({\bf F}_3)$&1&
${1\over 3}\epsilon_m(Q_8)+{1\over 6}\epsilon_m(C_2)$\\[1.5pt]
$SL_2({\bf F}_3)$&2&
${1\over 3}\epsilon_m(Q_8)+{2\over 3}\epsilon_m(C_2)$\\[1.5pt]
$SL_2({\bf F}_3)$&4&
${1\over 3}\epsilon_m(Q_8)+{5\over 3}\epsilon_m(C_2)$\\[1.5pt]
$SL_2({\bf F}_3)$&5&
${5\over 3}\epsilon_m(Q_8)+{5\over 6}\epsilon_m(C_2)$\\[1.5pt]
\hline
\end{tabular}
\end{center}
\end{minipage}
\begin{minipage}{2.25in}
\caption
{Values for $\delta_m(3)$.\label{tbl:wild3}}
\begin{center}
\vspace{-2ex}
\begin{tabular}
{|@{\hskip5pt}c@{\hskip5pt}|*{1}{@{\hskip5pt}c@{\hskip5pt}}|
@{\hskip5pt}c@{\hskip5pt}|}\hline
$\Phi_3$&$m=1$&formula\\\hline
$C_3,C_6$&2&$\epsilon_m(C_3)$\\[0.5pt]
$C_3\rtimes C_4$&1&${1\over 2}\epsilon_m(C_3)$\\[1.5pt]
$C_3\rtimes C_4$&3&${3\over 2}\epsilon_m(C_3)$\\[1.5pt]
\hline
\end{tabular}
\vspace{2.5ex}
\caption
{Values of $\Phi_2$.\label{tbl:phi2}}
\vspace{-2ex}
\begin{tabular}
{|@{\hskip5pt}c@{\hskip5pt}|
@{\hskip5pt}c@{\hskip5pt}|}\hline
$v_2(M)$&$\Phi_2$\\\hline

0&$C_1$ if $v_2(N)=0$ else $C_2$\\
2&$C_3$ if $v_2(N)=2$ else $C_6$\\
3,7&$SL_2({\bf F}_3)$\\
5&$Q_8$\\
8&$Q_8$ if $2^9|c_6(F)$ else $C_4$\\
\hline
\end{tabular}
\end{center}
\end{minipage}
\vspace{-8.0ex}
\end{table}

\subsection{The case of complex multiplication\label{sec:cm}}
\vspace{-0.5ex}
When $E$ has complex multiplication by an order of
some imaginary quadratic field~$K$, the situation simplifies since
we have \hbox{$L(E,s)=L(\psi,s-1/2)$}
for some\footnote{
 This is defined on ideals coprime
 to the conductor by~$\psi(z)=\chi(|z|)(z/|z|)$
 where $z$ is the primary generator of the ideal
 and $\chi$ is generally a quadratic Dirichlet character,
 but possibly cubic or sextic if $K=\Q\bigl(\sqrt{-3}\bigl)$,
 or quartic if~$K=\Q\bigl(\sqrt{-1}\bigr)$.
 When taking powers, we take $\chi^j$ to be the primitive Dirichlet character
 which induces~$\chi^j$.}
Hecke Gr\"ossencharacter~$\psi$.
For the symmetric powers we have the factorisation
\begin{equation}\label{cmfact}
L(\Sym^m E,s)=\prod_{i=0}^{m/2} L(\psi^{m-2i},s-m/2),
\end{equation}
where $\psi^0$ is the $\zeta$-function when $4|m$,
and when $2\|m$ it is $L(\theta_K,s)$ for the
quadratic character $\theta_K$ of the field~$K$.
Note that the local conductors and Euler factors
for each $L(\psi^j,s)$ can be computed iteratively
from \eqref{cmfact} since this information is known
for the left side from the previous subsections.
This factorisation reduces the computational complexity significantly,
as the individual conductors will be smaller than their product;
however, since there are more theoretical results in this case,
the data obtained will often lack novelty.
The factorisation~\eqref{cmfact} also implies that $L_{2u-1}(E,s)$
should vanish to high degree at~$s=u$, since each term
has about a 50\% chance of having odd functional equation.
We found some examples where $L(\psi^3,s)$, $L(\psi^5,s)$,
or $L(\psi^7,s)$ has a double zero at the central point,
but we know of no such triple zeros.

\section{Global considerations and computational techniques}
We now give our method for computing special values of the symmetric
power \hbox{$L$-functions} defined above.
To do this, we complete the \hbox{$L$-function} with
a \hbox{$\Gamma$-factor} corresponding to the prime at infinity, and then
use the (conjectural) functional equation in conjunction with the
method of Lavrik~\cite{lavrik} to write the special value as
a ``rapidly-converging'' series whose summands involve
inverse Mellin transforms related to the \hbox{$\Gamma$-factor.}
First we digress on poles of our~\hbox{$L$-functions}.

\subsection{Poles of $L$-functions}
It is conjectured that $L_m(E,s)$ has an entire continuation,
except when $4|m$ and $E$ has complex multiplication (CM) there is
a pole at $s=1+m/2$, which is the edge of the critical strip.\footnote
{The case of $m=4$ follows as a corollary of work
 of Kim~\cite[Corollary~7.3.4]{kim}.}
We give an explanation of this expectation
from the standpoint of analytic number theory;
it is likely that a different argument could be given
via representation theory.
We write each Euler factor as~$U_m(p;s)=\bigl(1-b_m(p)/p^s+\cdots\bigr)^{-1}$
and as $s\rightarrow 1+m/2$ we have $\log L_m(s)\sim \sum_p b_m(p)/p^s$.
We will now compute that the conjectural Sato-Tate distribution~\cite{tate}
implies that the average value of $b_m(p)$ is~0,
while for CM curves the Hecke distribution~\cite{hecke}
will yield an average value for $b_m(p)$ of $p^{m/2}$ when~$4|m$.

Similar to~\eqref{eqn:sigma}, for a good prime $p$ we have
$b_m(p)=\sum_{i=0}^{m/2} {m-i\choose i} a_p^{m-2i}(-p)^i$.
The Sato-Tate and Hecke distributions imply that the average values of
the $k$th power of $a_p$ are given by
$$\langle a_p^k\rangle=(2\sqrt p)^k
{\int_0^\pi (\cos\theta)^k\,(\sin\theta)^2\,d\theta\over
\int_0^\pi (\sin\theta)^2\,d\theta}
\>\>\>\>{\rm and}\>\>\>\>
\langle a_p^k\rangle_{CM}=(2\sqrt p)^k
{\int_0^\pi (\cos\theta)^k\,d\theta\over 2\int_0^\pi d\theta}.$$
We have $\langle a_p^k\rangle=0$ for $k$~odd;
for even~$k$ the Wallis formula \cite{wallis} implies
$$\int_0^\pi (\cos\theta)^k(\sin\theta)^2\,d\theta=
{\pi (k-1)!!\over k!!}-{\pi (k+1)!!\over(k+2)!!}={\pi (k-1)!!\over (k+2)!!},$$
so that $\langle a_p^k\rangle$ is $(2\sqrt p)^k{2(k-1)!!\over (k+2)!!}$.
An induction exercise shows that this implies $\langle b_m(p)\rangle=0$
when $E$ does not have~CM.
We also have $\langle a_p^k\rangle_{CM}=(2\sqrt p)^k{(k-1)!!\over 2\cdot k!!}$
for even~$k$, and again an inductive calculation shows
that $\langle b_m(p)\rangle_{CM}=p^{m/2}$ when $4|m$ and is zero otherwise.
This behaviour immediately implies the aforementioned
conjecture about the poles of~$L_m(E,s)$ at~$s=1+m/2$.

\subsection{Global considerations}
Let $\Lambda_m(E,s)=C_m^s\gamma_m(s)L_m(E,s)$,
where $C_m^2=N_m/(2\pi)^{m+1}$ for $m$ odd and is twice this for $m$ even.
For $m$ odd we write $m=2u-1$ and for $m$ even we write $m=2v$;
then from \cite[\S 5.3]{deligne} we have respectively either
$$\gamma_m(s)=\prod_{i=0}^{u-1} \Gamma(s-i)
\qquad{\rm or}\qquad
\gamma_m(s)=
\Gamma\bigl(s/2-\lfloor v/2\rfloor\bigr)\prod_{i=0}^{v-1} \Gamma(s-i).$$
When $4|m$ and $E$ has CM, we multiply $\gamma_m(s)$ by~$(s-v)(s-v-1)$.
We expect $\Lambda_m(E,s)$ to have an entire continuation and satisfy a
functional equation $\Lambda_m(E,s)=w_m \Lambda_m(E,m+1-s)$
for some $w_m=\pm 1$.
The works of Kim and Shahidi \cite{ks2} establish
parts of this conjecture.\footnote{
 The full conjecture\footnotemark
 follows from Langlands functoriality~\cite{langlands}.
 In the CM case, the functional equation follows
 from the factorisation~\eqref{cmfact} and the work of Hecke~\cite{hecke}.}
\footnotetext{Added in proof: A recent preprint \cite{taylor} on Taylor's
 webpage shows the meromorphic continuation and functional equation
 for all symmetric powers for curves with~$j\not\in\Z$.\looseness=-1}
We can find $w_m$ via experiment as described in
Section~\ref{section:computations},
but we can also try to determine $w_m$ theoretically.

\vspace*{-4pt}
\subsection{Digression on local root numbers}
\vspace*{-2pt}
The sign $w_m$ can theoretically be determined via local computations
as in~\cite{deligne:rootnumber}, but this is non-trivial to implement
algorithmically, especially when~$p=2,3$.
We expect to have a factorisation $w_m=\prod_p w_m(p)$
where the product is over bad primes~$p$ including infinity.
For $m$ even, the very general work of
Saito \cite{saito} can then be used to show\footnote{
 The work of Fr\"ohlich and Queyrut~\cite{frohlich-queyrut}
 and Deligne~\cite{deligne:invent} might give a direct argument.}
that~$w_m=+1$, so we assume that $m$ is odd.
From~\hbox{\cite[\S 5.3]{deligne}}
we have $w_m({\infty})=-\bigl({-2\over m}\bigr)$;
combined with the relation $w_m(p)=w_1(p)^m$ for primes $p$ of
multiplicative reduction, this gives the right sign for semistable curves.
The potentially multiplicative case has~$w_m(p)=w_1(p)^{(m+1)/2}$.

In the additive cases,
the first author \cite{martin} has used the work of Rohrlich \cite{rohrlich}
to compute the sign for $p\ge 5$. We get that\footnote{
 Since we are assuming that $m$ is odd,
 the exponent is just $(m+1)/2$ unless~$\Phi_p=C_3$.}
$w_m(p)=w_1(p)^{\epsilon_m(\Phi_p)/2}$,
and $w_1(p)$ is listed in~\cite{rohrlich}.
For $p=2,3$ the value of $w_1(p)$ is given\footnote{
 Note the third case in Table~1 of \cite{halberstadt} needs a
 \textit{Condition sp\'eciale} of~$c_4'\equiv 3\pmod{4}$.}
by Halberstadt~\cite{halberstadt},
and our experiments for higher (odd) powers indicate that
$w_m(2)=\eta_2 w_1(2)^{\epsilon_m(\Phi_p)/2}$
where $\eta_2=-1$ if $v_2(N)$ is odd and $m\equiv 3\pmod{8}$
and else~$\eta_2=+1$, while the expected values of
$w_m(3)$ are given in Table~\ref{tbl:wm3}.

\vspace{-3ex}
\begin{table}[h]
\caption
{Experimental values for $w_m(3)$ (periodic mod~12 in~$m$.\label{tbl:wm3})}
\begin{center}
\vspace{-2ex}
\begin{tabular}
{|@{\hskip4pt}c@{\hskip4pt}|*{6}{@{\hskip3pt}c@{\hskip3pt}}
|@{\hskip4pt}c@{\hskip4pt}|*{6}{@{\hskip3pt}c@{\hskip3pt}}
|@{\hskip4pt}c@{\hskip4pt}|*{6}{@{\hskip3pt}c@{\hskip3pt}}|}\hline
$\Phi_3$&1&3&5&7&9&11&
$\Phi_3$&1&3&5&7&9&11&
$\Phi_3$&1&3&5&7&9&11\\\hline

$C_3,C_4$&$+$&$+$&$+$&$+$&$+$&$+$&
$C_6$&$+$&$-$&$-$&$-$&$+$&$+$&
$C_3\rtimes C_4$&$+$&$+$&$-$&$+$&$+$&$+$\\
$C_2$&$-$&$+$&$-$&$+$&$-$&$+$&
$C_6$&$-$&$+$&$-$&$+$&$-$&$+$&
$C_3\rtimes C_4$&$-$&$-$&$-$&$-$&$-$&$+$\\
\hline
\end{tabular}
\end{center}
\vspace{-7.5ex}
\end{table}

\vspace*{-4pt}
\subsection{Computations\label{section:computations}}
\vspace*{-2pt}
From~\cite{lavrik},~\hbox{\cite[Appendix B]{cohen},}~or~\cite{dokchitser},
the assumption of the functional equation
$\Lambda_m(E,s)=w_m\Lambda_m(E,m+1-s)$
allows us to compute (to a given precision) any value/derivative
$\Lambda_m^{(d)}(E,s)$ in time proportional to~$C_m\approx \sqrt{N_m}$.
Additionally, numerical tests on the functional equation
arise naturally from the method.

We follow \hbox{\cite[\S~7,~p.~119ff]{BST}.}
Suppose we have $\Lambda_m(s)=w_m\Lambda(m+1-s)$,
and the $d$th derivative is the first one that is nonzero at~$s=\kappa$.
Our main interest is in $\kappa=u$ for $m=2u-1$ and $\kappa=v+1$ for $m=2v$,
and we note that $d=0$ for even~$m$.
Via Cauchy's residue theorem, for every real $A>0$ we have
$${\Lambda_m^{(d)}(\kappa)\over d!}=
{1\over 2\pi i}\biggl(\int_{(\delta)}-\int_{(-\delta)}\biggr)
{\Lambda_m(z+\kappa)\over z^{d+1}}{dz\over A^z},$$
where $\delta$ is small and positive and $\int_{(\sigma)}$
is the integral along~$\Re z=\sigma$.
In the second integral we change variables $z\rightarrow -z$
and apply the functional equation.
Then we write $\kappa+\lambda=m+1$, move both contours sufficiently
far to right (say~$\Re z=2m$) and expand $\Lambda_m$
in terms of the $L$-function to get
\begin{align*}
{\Lambda_m^{(d)}(\kappa)\over d!}=
\int_{(2m)}C_m^{z+\kappa}&\gamma_m(z+\kappa)
\sum_{n=1}^\infty {b_m(n)\over n^{z+\kappa}}
{1\over z^{d+1}}{dz\over 2\pi i\,A^z}\\
&+(-1)^dw_m\int_{(2m)}C_m^{z+\lambda}\gamma_m(z+\lambda)
\sum_{n=1}^\infty {b_m(n)\over n^{z+\lambda}}
{1\over z^{d+1}}{A^z\,dz\over 2\pi i}.
\end{align*}
Thus we get that
$${\Lambda_m^{(d)}(\kappa)\over d!}=
C_m^{\kappa}\sum_{n=1}^\infty {b_m(n)\over n^\kappa}
F^d_m\bigl(\kappa;{\textstyle{n\over AC_m}}\bigr)
+(-1)^dw_m C_m^{\lambda}\sum_{n=1}^\infty {b_m(n)\over n^{\lambda}}
F^d_m\bigl(\lambda;{\textstyle{nA\over C_m}}\bigr),$$
where
$$F_m^d(\mu;x)=
\int_{(2m)} {\gamma_m(z+\mu)\over z^{d+1} x^z}{dz\over 2\pi i}.$$
The $F_m^d(\mu;x)$-functions are ``rapidly decreasing''
inverse Mellin transforms.
Note that we have
$L_m^{(d)}(\kappa)=\Lambda_m^{(d)}(\kappa)/\gamma_m(\kappa)C_m^\kappa$,
and so can recover the $L$-value as desired.
The parameter $A$ allows us to test the functional equation;
if we compute $\Lambda_m^{(d)}(\kappa)$ to a given precision for $A=1$
and $A=9/8$, we expect disparate answers
if we have the wrong Euler factors or sign~$w_m$.

We compute $F_m^d(\mu;x)$ as a sum of residues at poles in the left half-plane,
the first pole being at~$z=0$, following~\cite{delaunay}.
We need to calculate Laurent series expansions
of the $\Gamma$-factors about the poles.\footnote{
 When $\Lambda_m(E,s)$ has a pole
 the factor $\gamma_m(s)$ has two additional linear factors
 (which are easily handled). But in this case it is better
 to use the factorisation \eqref{cmfact}.}
We let $\zeta(1)$ denote Euler's constant $\gamma\approx 0.577$,
and define $H_1(n)=1$ for all~$n$, and $H_k(1)=\sum_{i=1}^k 1/i$ for all~$k$,
and recursively define $H_k(n)=H_{k-1}(n)+H_{k}(n-1)/k$ for~$n,k\ge 2$.
At a pole $z=-k$ for $k$ a nonnegative integer, we have the Laurent expansion
$$\Gamma(z)=
{(-1)^k\over k!(z+k)}\biggl(1+\sum_{n=1}^\infty H_k(n)(z+k)^n\biggr)
\exp\biggl(\sum_{n=1}^\infty {(-1)^n\zeta(n)\over n}(z+k)^n\biggr),$$
and for $k$ a negative integer (these only occur for a few cases)
we can use the relation $z\Gamma(z)=\Gamma(z+1)$ to shift.
To expand $\Gamma(z/2)$ around an odd integer~\hbox{$z=-k$,}
we use the duplication formula
$\Gamma(z)=\Gamma(z/2)\Gamma\bigl({z+1\over 2}\bigr){\sqrt\pi\over 2^{z-1}}$
to replace $\Gamma(z/2)$ by a quotient of $\Gamma$-factors
that can each be expanded as above.
The trick works in reverse to expand $\Gamma\bigl({z+1\over 2}\bigr)$
about an even integer~\hbox{$z=-k$.}
We also have the series expansions for $2^z$ and $1/z$ about~$z=-k$ given by
$$2^z=2^{-k}\sum_{n=0}^\infty {(\log 2)^n\over n!}(z+k)^n
\quad{\rm and}\quad
{1\over z}=-{1\over k}-\sum_{n=1}^\infty {(z+k)^n\over k^{n+1}}
\>\>(\text{for } k\neq 0).$$

Since these $F_m^d(\mu;x)$ functions are (except for CM)
independent of the curve,
we pre-computed a large mesh of values and derivatives of these functions,
and then in our programme we compute via local power series.
Thus, unlike the setting of Dokchitser~\cite{dokchitser},
we are not worried too much about the cost of computing $F_m^d(\mu;x)$
for large~$x$ via a massively-cancelling series expansion,
since we only do this in our pre-computations.
For each implemented function we have its value
and first 35 derivatives for all $x=i2^k/32$ for $32\le i\le 63$
for $k$ in some range, such as $-3\le k\le 19$.
For sufficiently small~$x$ we just use the log-power-series expansion.
The choice of 35 derivatives combined with the maximal radius
of $x/64$ for expansions about~$x$ implies that our
maximal precision is around $35\times 6=210$ bits.
When working to a lower precision, we need not sum so many terms
in the local power series. Note that $F_m^d(\mu;x)$ dies off roughly
like $\exp(-x^{2/(m+1)})$, and thus it is difficult to do
high precision calculations for $m$ large.

\setbox4\hbox{\verb+www.maths.bris.ac.uk/~mamjw+}
To compute the meshes of inverse Mellin transforms described above,
we used PARI/GP~\cite{pari}, which can compute to arbitrary precision.
However, PARI/GP was too slow to use when actually computing the $L$-values;
instead we used a C-based adaption of Bailey's quad-double
package~\cite{bailey-qd}, which provides up to 212~bits of precision
while remaining fairly fast.\footnote{
 The SYMPOW package can be obtained from \box4}

\section{Results}
We tested the functional equation (via the above method of comparing
the computed values for $A=1$ and $A=9/8$)
for odd symmetric powers $m=2u-1$ at the central point~$\kappa=u$,
and for even symmetric powers $m=2v$
at the edge of the critical strip~$\kappa=v+1$.
We did this for all non-CM isogeny classes in Cremona's
database~\cite{cremona} with conductor less than~130000;
this took about 3 months on a cluster of 48 computers
(each running at about 1 Ghz).

We computed as many as $10^8$ terms of the
various \hbox{$L$-series} for each curve,
which was always sufficient to check the functional equation
of the third symmetric power to about six decimal digits.\footnote{
 In about 0.3\% of the cases, the computations for both the zeroth
 and first derivatives showed no discrepancy with $A=1$ and $A=9/8$;
 this coincidence is to be expected on probabilistic grounds, and for these
 cases we computed to higher precision to get an experimental
 confirmation of the sign of the functional equation.}
In all cases, we found the expected functional equation
to hold to the precision of our calculation.
The results for the order of vanishing (at the central point) for odd powers
appear in the left half of Table~\ref{table:results}.
The right half lists how many tests\footnote{
 We need not compute even powers when there is
 a lack of quadratic-twist-minimality.}
we did for other\footnote{
 We did not test the fourth symmetric power,
 as the work of Kim \cite{kim} proves the validity of the
 functional equation in this case. Since there is
 no critical value, a calculation would do little more
 than verify that our claimed Euler factors are correct.}
symmetric powers (again to six digits of precision).
There are less data for higher symmetric powers due to
our imposed limit of $10^8$ terms in the $L$-series computations,
but since the symmetric power conductors for curves with exotic
inertia groups often do not grow so rapidly, we can still test
quite high powers in some cases.

\eject
\begin{table}[h]
\caption
{\small{Test-counts (right) and data for order
of vanishing (non-CM isogeny classes)}\looseness=-1\label{table:results}}
\vspace{-1.5ex}
\hfil
\begin{minipage}{2.55in}
\begin{center}
\begin{tabular*}{2.55in}
{@{\extracolsep{\fill}}|c@{\,\,}|r@{\,\,}|
r@{\,\,}|r@{\,\,}|r@{\,\,}|r@{\,\,}|r@{\,\,}|}\hline
\,{$m$}&Tested&Order $0$&Order $1$&$2$\ \ \ &$3$\ \ &$4$\, \\\hline
\,$1$&567735&216912&288128&61787&908&0\\
\,$3$&567735&262751&287281&16782&905&16\\
\,$5$&46105&22448&23076&569&12&0\\
\,$7$&3573&1931&1616&25&1&0\\
\,$9$&947&542&400&5&0&0\\
\,$11$&134&51&82&1&0&0\\\hline
\end{tabular*}
\end{center}
\end{minipage}
\hfil
\begin{minipage}{0.7in}
\begin{center}
\begin{tabular*}{0.7in}{@{\extracolsep{\fill}}|c@{\,\,}|r@{\,\,}|}\hline
\,$m$&\# tests\\\hline
\,$6$&4953\\
\,$8$&1259\\
\,$10$&190\\
\,$12$&142\\
\,$13$&5 even\\
\,$13$&30 odd\\\hline
\end{tabular*}
\end{center}
\end{minipage}
\hfil
\begin{minipage}{0.7in}
\begin{center}
\begin{tabular*}{0.7in}{@{\extracolsep{\fill}}|c@{\,\,}|r@{\,\,}|}\hline
\,$m$&\# tests\\\hline
\,$14$&26\\
\,$15$&1 even\\
\,$15$&16 odd\\
\,$16$&8\\
\,$17$&3 odd\\
\,$18$&2\\\hline
\end{tabular*}
\end{center}
\end{minipage}
\end{table}
\vspace{-4ex}

Buhler, Schoen, and Top \cite{BST} already listed
2379b1 and 31605ba1 as 2 examples
of (suspected) 4th order zeros for the symmetric cube.
We found 14 more, but no examples of 5th order zeros.
For higher powers, we found examples of 3rd order zeros
for the 5th and 7th powers, and 2nd order zeros for the
9th, 11th, and 13th powers, though as noted above,
we cannot obtain as much data for higher powers.\footnote
{Given that we only computed the $L$-value of the 13th symmetric power
 for five curves of even sign,
 to find one that has a double-order zero is rather surprising.
 Higher-order zeros were checked to 12 digits;
 the smallest ``nonzero'' value was~$\approx 2.9\cdot 10^{-8}$.}
We list the Cremona labels for the
isogeny classes in Table~\ref{table:vanishing}.

\vspace{-3ex}
\begin{table}[h]
\caption{\small{Experimentally observed high order vanishings
(non-CM isogeny classes)}\looseness=-1\label{table:vanishing}}
\vspace{-4ex}%
\begin{center}
\begin{tabular}{|@{\hskip3pt}c@{\hskip3pt}|@{\hskip3pt}l@{\hskip3pt}|}\hline
ord&format is {\bf power}:label(s)\\\hline
4th& {\bf 3}:2379b 5423a 10336d 29862s 31605ba 37352d 46035a 48807b 55053a \\
&{\bf 3}:59885g 64728a 82215d 91827a 97448a 104160bm 115830a \\
3rd& {\bf 5}:816b 2340i 2432d 3776h 5248a 6480t 7950w 8640bl
16698s 16848r \\
& {\bf 5}:18816n 57024du\hskip18pt{\bf 7}:176a \hfil\\
2nd& {\bf 7}:128b 160a 192a 198b 200e
320b 360b 425a 576b 726g 756b 1440a \\
& {\bf 7}:1568i 1600b 2304g 3267f
3600h 3600j 3600n 3888e 4225m 6272d \\
& {\bf 7}:11552r 15876f 21168g \hfil {\bf 9}:40a 96a 162b 324d 338b
\hfil {\bf 11}:162b \hfil {\bf 13}:324c\\\hline
\end{tabular}
\end{center}
\end{table}
\vspace{-6ex}

We also looked at extra vanishings of the 3rd symmetric power
in a quadratic twist family.
We took $E$ as 11a3:$[0,-1,1,0,0]$ and computed the twisted
central value $L_3(E_d,2)$ or central derivative $L'_3(E_d,2)$
for fundamental discriminants~$|d|<5000$.
We found 58 double zeros (to 9 digits) and one triple zero ($d=3720$).
A larger experiment (for $|d|<10^5$) for 10 different CM curves
found (proportionately) fewer double zeros and no triple zeros.

Finally, we used higher-precision calculations to obtain
the Bloch-Kato numbers of equation \eqref{eqn:special}
for various symmetric powers of some non-CM curves of small conductor
(see Table~\ref{table:BK}). More on the arithmetic significance
of these quotients will appear elsewhere. In some cases, we were able
to lessen the precision because it was known that a large power
of a small prime divided the numerator.

\begin{table}[ht]
\caption%
{Selected Bloch-Kato numbers for various powers and curves\label{table:BK}}
\vspace{-2ex}
\hfil
\begin{minipage}{90pt}
\begin{center}
\begin{tabular}{|@{\hskip4pt}l@{\hskip4pt}l@{\hskip4pt}|}\hline
\rlap{\centerline{\textbf{5th powers}}} & \\
20a2& $2^9$ \\
37a1& $2^9$ \\
43a1& $2^7 5$\\
44a1& $2^{17}$\\
\hline
\rlap{\centerline{\textbf{6th powers}}} & \\
11a3& $2^4 5$ \\
14a4& $2^9 3$ \\
15a8& $2^{10}$ \\
17a4& $2^{12}$ \\
19a3& $2^4 3^3 5^2$ \\
20a2& $2^{17}/3$ \\
24a4& $2^{17}/3$ \\
26a3& $2^{7} 3\cdot 5\cdot 23$ \\
26b1& $2^{7} 3\cdot 7^3\cdot 23$ \\
30a1& $2^{15} 3^3 7$ \\
33a2& $2^{17} 3\cdot 5\cdot 7$ \\
34a1& $2^{13} 3^3 59$ \\
35a3& $2^8 3\cdot  7^2 31$ \\
37a1& $2^9 3^4 7$ \\
37b3& $2^7 3^4 467$ \\
38a3& $2^7 3^4 5\cdot 11\cdot 137$ \\
38b1& $2^7 3^2 5^2 13\cdot 31$ \\
39a1& $2^{20} 3^2 7$ \\
40a3& $2^{20} 7/3$ \\
42a1& $2^{19} 3^2 7\cdot 19$\\
43a1& $2^6 3^2 1697$ \\
44a1& $2^{21} 5\cdot 31/3$ \\
46a1& $2^9 5\cdot 23\cdot 30661$ \\
50a1& $2^3 5^{11} / 3$ \\
51a1& $2^9 3^3 4517$ \\
\hline
\end{tabular}
\end{center}
\end{minipage}
\hfil\hskip10pt%
\begin{minipage}{95pt}
\begin{center}
\begin{tabular}{|@{\hskip4pt}l@{\hskip4pt}l@{\hskip4pt}|}\hline
\rlap{\centerline{\textbf{7th powers}}} & \\
24a4& $2^{23} 7/3$ \\
37a1& $2^{13} 3\cdot 5$ \\
43a1& $2^{17} 3\cdot 5$ \\
\hline
\rlap{\centerline{\textbf{9th powers}}} & \\
11a3& $2^{12}$ \\
14a4& $2^{14} 3^4 5$ \\
15a8& $2^{16}$ \\
17a4& $2^{16} 3^6 5$ \\
19a3& $2^{19} 3^2 5$ \\
21a4& $2^{20} 5\cdot 59^2$ \\
24a4& $2^{38}/9$ \\
26a3& $2^{11} 3^4 5\cdot 7^4$ \\
26b1& $2^{11} 3^2 5\cdot 7^3 1933^2$ \\
30a1& $2^{16} 3^5 5^5 37^2$ \\
33a2& $2^{24} 5\cdot 107^2 167^2$ \\
34a1& $2^{23} 3^5 5\cdot 7^2 53^2$\\
35a3& $2^{25} 3^4 5$ \\
37b3& $2^{20} 3^2 5\cdot 7^2 53^2$ \\
38a3& $2^{11} 3^{14} 5\cdot 19^2$ \\
38b1& $2^{11} 5^8 109^2$ \\
39a1& $2^{40} 3^2 5 7^4$ \\
40a3& $0$ \\
42a1& $2^{25} 5\cdot 223^2 241^2$ \\
44a1& $2^{47} 3$ \\
45a1& $2^{16} 3^{19} 5^2 7\cdot 13^2$ \\
46a1& $2^{14} 3^{10} 5^3 14071^2$\\
48a4& $2^{43} 5$ \\
50a1& $2^5 3\cdot 5^{22}$ \\
54a3& $2^9 3^{24}$ \\
54b1& $2^7 3^{25} 5$ \\
\hline
\end{tabular}
\end{center}
\end{minipage}
\hfil\hskip10pt%
\begin{minipage}{125pt}
\begin{center}
\begin{tabular}{|@{\hskip4pt}l@{\hskip4pt}l@{\hskip4pt}|}\hline
\rlap{\centerline{\textbf{10th powers}}} & \\
11a3& $2^{14} 5\cdot 22453/3$ \\
14a4& $2^{16} 3^3 5\cdot 6691$ \\
15a8& $2^{26} 5\cdot 541$ \\
17a4& $2^{23} 3^2 7\cdot 11\cdot 227$ \\
19a3& $2^{14} 3^2 47\cdot 179\cdot 5023$ \\
20a2& $2^{44} 53/3$ \\
21a4& $2^{28} 3^7 5^2 29$ \\
24a4& $2^{49} 13/9$ \\
26a3& $2^{19} 3^5 7\cdot 47\cdot 1787$ \\
26b1& $2^{19} 3^3 5^2 7^3 127\cdot 2102831$ \\
40a3& $2^{54} 5\cdot 683$ \\
44a1& $2^{56} 5\cdot 11\cdot 215447 / 3$ \\
50a1& $2^{11} 5^{28} 7/3$ \\
52a2& $2^{44} 3^3 5\cdot 7\cdot 19\cdot 279751$\\
54b1& $2^{14} 3^{35} 7$ \\
56a1& $2^{66} 3^2 5\cdot 11\cdot 71$\\
75c1& $2^{14} 5^{28}\cdot 31\cdot 41\cdot 61/9$ \\
96b1& $2^{84} 197/3$ \\
99a1& $2^{18} 3^{31} 5^3 7\cdot 1367$\\
\hline
\rlap{\centerline{\textbf{11th powers}}} & \\
11a3& $2^{26} 5^4/3$ \\
14a4& $2^{23} 3^5 5^2 7^2 11^2$ \\
15a8& $2^{29} 11^2 23^2/3$ \\
17a4& $2^{26} 3^{11}$ \\
21a4& $2^{36} 11^2 211^2 / 3 $ \\
24a4& $2^{57} 13^2/45$ \\
48a4& $2^{70} 11^2/3$ \\
54b1& $2^{20} 3^{41}$\\
56a1& $2^{74} 11^2 / 5 $ \\
72a1& $2^{58} 3^{28} 59^2/5$\\
\hline
\end{tabular}
\end{center}
\end{minipage}
\vspace{-4ex}
\end{table}

\vspace{-2ex}
\subsection{Other directions}
In this work, we looked at symmetric powers for weight~2 modular forms.
Delaunay has done some computations \cite{delaunay}
for modular forms of higher weight;
in that case, the work of Deligne again tells us where
to expect critical values, and the experiments confirm that we do indeed
get small-denominator rationals after proper normalisation.
We looked at critical values at the edge and center of the
critical strip, whereas we expect $L$-functions evaluated at
other integers to take special values related to $K$-theory;
see
\cite{bloch-grayson,mestre-schappacher,dokchitser-dejeu-zagier,zagier-gangl}
for examples.
The programmes written for this paper
are readily modifiable to compute other special values.
The main advantage that our methods have over those of
Dokchitser \cite{dokchitser} is that we fixed the
$\Gamma$-factors and the $L$-values of interest, which
then allowed a large pre-computation for the inverse
Mellin transforms; if we wanted (say) to compute
zeros of $L$-functions (as~with~\cite{rubinstein}),
our method would not be as useful.

Finally, the thesis of Booker \cite{booker} takes another approach
to some of the questions we considered. The scope is much more broad,
as it considers not only numerical tests of modularity, but also tests
of GRH~(\S 3.4), recovery of unknown Euler factors possibly using
twists~(\S 5.1), and also high symmetric powers~(\S 7.2).

\end{document}